\documentclass[11pt]{amsart}

\usepackage{amsmath,amssymb,amsfonts}

\usepackage{caption2}
\usepackage[]{hyperref}

\usepackage[pdftex]{color,graphicx}
\usepackage{multicol}
\usepackage{verbatim}

\numberwithin{equation}{section}

\setcaptionmargin{1cm}

\newcommand{\be}{\begin{equation}}
\newcommand{\ee}{\end{equation}}

\newcommand{\sign}{\text{sign}}

\usepackage{color}
\definecolor{rosso}{cmyk}{0,1,1,0.4}
\definecolor{rossos}{cmyk}{0,1,1,0.55}
\definecolor{rossoc}{cmyk}{0,1,1,0.2}

\begin{document}

\begin{center}
{
 \bf

Phase variation  and angular momentum of  the Riemann  and Dirichlet $\xi$  functions. %
  }  
\vskip 0.2cm {\large

Giovanni  Lodone \footnote{
giolodone3@gmail.com}  %
.} 
\\[0.2cm]

{\bf

ABSTRACT    } 
\\[0.2cm]

\vskip 0.2cm {\large

The concept of angular momentum  
 is used to find new RH equivalence statements, and, generalize some known results  from Riemann to Dirichlet primitive $\xi$ functions } 
\\[0.2cm]

\vskip 0.1cm

\end{center}

{

MSC-Class  :   11M06, 11M26, 11M99

  } 
\vspace{0.1cm}



     \tableofcontents

\section{ Introduction}

In \cite {Giovanni Lodone Nov2024}  and \cite {Giovanni Lodone Dec2024}  it is used a function called, angular momentum (\ref{AngMomDef}),  in order to not reinvent existent named quantities, to attempt an RH proof for two types of Dirichlet function: Riemann $\zeta(s)$ function, and, $L(s,\chi_{odd \ primitive})$ Dirichlet functions. See (\cite  {Davenport1980},  \cite {Apostol:1976}, and,   \cite [p.~1] {Peter Sarnak}  ). 

\noindent In this paper %
we show that this concept is useful also to give an unitary treatment  of  some aspects of Dirichlet functions. Besides some known results about Riemann $\xi$ function can be extended easily to primitive Dirichlet function. We apologize for unavoidable overlapping with previous papers.


\noindent We express the complex argument $s$, as usual, as:
\be \label{1p2}
 s=\frac{ 1}{2}+\epsilon +it 
 \ee
where the real values $\epsilon$ and $t$ are respectively the distance from the critical line ( i.e $\Re(s) = 1/2$) and the imaginary part of $s$, while,  $0 < \Re(s)  < 1$ is the critical strip.

  \noindent As well known Rieman $\zeta(s)$ function, \cite{Riemann:1859ar} (see also  \cite{Bombieri:2000cz}, \cite{Borwein:2007cz}, 
\cite[p.~229] {Edwards:1974cz}), is defined as the analytic continuation of the function:

\be \label {1p1}
\zeta(s)=\sum _{n=1}^\infty \frac{1}{n^s} = \prod_{\forall p} \frac{1}{1-\frac{1}{p^s}} \quad \mbox{ ($p$ prime)}  \quad \quad \Re(s)>1
\ee

\noindent { \bf [RH ( Riemann Hypothesis) ]$\equiv$[  the analytic continuation of (\ref {1p1})  has all  non trivial zeros,  $s_0$,
  on critical-line]}.

\noindent Using   Dirichlet characters $\chi(n)$, L functions  (\cite  {Davenport1980} and  \cite {Apostol:1976},   \cite [p.~1] {Peter Sarnak}  ) are defined as :
   
   \be \label {DirichletLFunc}
\L(s.\chi)=\sum _{n=1}^\infty  \   \chi(n)/n^s = \prod_{\forall p} \frac{1}{1-\frac{\chi(p)}{p^s}} \quad p \quad prime  \quad  ,  \quad  \quad \quad \Re(s)>1
\ee

{ \bf [GRH (Generalized Riemann Hypothesis)]$\equiv$[
all the zeros inside critical strip of the (\ref{DirichletLFunc}), or of their analytic continuation,  are 
  on critical-line]}, i.e $ \forall s_0 $ inside critical strip $: \Re(s_0) = 1/2$ for (\ref {1p1}) and (\ref{DirichletLFunc}).

\noindent Because of (\ref{FuncEquForXiL}) we focus  mainly on $\chi_{primitive}$, defined in  \cite [p.~168]{Apostol:1976}.  
Companion functions   $\xi(s)$,  (  \cite [p.~62] {Davenport1980} ), and $\xi(s,\chi_{primitive} )$, (  \cite [p.~71] {Davenport1980} ), are defined :

\
     
     \be \label {XiForLfunc} 
\xi(s)= \Gamma \left(\frac{s}{2} +1\right)(s-1) \frac {\zeta(s)}{\pi^{s/2}} \ ;  \xi(s,\chi) = \left( \frac{q}{\pi}   \right)^{\frac{s+\alpha}{2}} \Gamma  \left( \frac{s+\alpha}{2} \right) L(s,\chi)  
 \ee

 \noindent  where $ q= congruences \quad modulus $,  with  $\alpha=0$ if $\chi_{primitive}(-1) =1 $, even character,  and,   $\alpha=1$ if $\chi_{primitive}(-1) =-1$, odd character. 
Both have the same zeros of the $\zeta(s)$ or $L(s,\chi)$ functions  respectively in the critical strip \cite{Riemann:1859ar} \cite{Davenport1980} \cite{Edwards:1974cz}.
     \noindent       Following    \cite[p.~14]{Borwein:2007cz}, and 
     \cite [p.~71] {Davenport1980} 
      the functional equation for $\xi(s)$ and %
       $\xi(s,\chi)$  
      can be  defined  :

        \be  \label {FuncEquForXiL}
       \quad \xi(1-s)=\xi(s)  \ \ \ ; \ \ \
        \xi(1-s,\bar {\chi})=\frac{i^\alpha q^{\frac{1}{2}}}{\tau(\chi)} \xi(s,\chi) \quad with \quad \chi_{primitive} \ \ %
         \quad
        \ee
          
           \noindent  Where $\bar {\chi}$ is the complex conjugate of $ \chi$ and
              \noindent $\tau(\chi)$ is the gauss sum $\tau(\chi) :=  \sum_{m=1}^k \chi(m) e^{2 \pi i m/k}$  
              see \cite [p.~65] {Davenport1980} ,           
        \noindent        or \cite [p.~165]  {Apostol:1976} . %
                
              \noindent  For \cite [p.~168]  {Apostol:1976}, and, \cite [p.~66] {Davenport1980}:  $%
              |\tau(\chi)  |^2= q$
                where $q$ is the arithmetic congruence modulus, (\ref{FuncEquForXiL})  leads to mirror symmetry around $\Re s =1/2$ and real axis for $\xi(s)$: %
\be \label {ConiugioDxSx}
\xi \left (\frac{ 1}{2}+\epsilon +it  \right) = \bar {   \xi  }\left (\frac{ 1}{2}-\epsilon +it  \right )  
\ee

   \noindent  Instead if  $L(s,\chi)$   (\ref{DirichletLFunc}) or  $\xi(s,\chi)$  have off-critical line zeros they are bound to respect  in general
   only mirror symmetry with respect to critical line,  i.e. a zero in $(t,\epsilon)$ means a zero  in   $(t,-\epsilon)$. 
        See  \cite[p.~83] {Davenport1980}.  %
    In  \cite[p.~80] {Davenport1980} and in\cite[p.~82] {Davenport1980}    similar expressions  are given for Hadamard infinite product  of %
     (\ref  {XiForLfunc}):

\be \label {InfiniteProdHad}
\xi(s) = e^{g(s)} \prod_{\rho \quad paired} \left(  1-\frac{s}{\rho} \right) e^{s/\rho}
\quad ;\quad 
\xi(s,\chi) = e^{g(s,\chi)} \prod_{|\Im \rho| <T\rightarrow \infty} \left(  1-\frac{s}{\rho} \right) e^{s/\rho} 
\ee
where $g(z,\chi)=A(\chi)+B(\chi) z$.

\be  \label {BAndRoots}
For \quad \xi(s) :\ \ B= - \sum_{m=-\infty}^ {m=+\infty} \frac {1}{\rho_m} \quad;\quad  For \quad   \xi(s,\chi):\ \ \Re[B(\chi)[= - \sum_{\Im \rho_m>-T}^ {\Im \rho_m<T}\Re\left[ \frac {1}{\rho_m} \right]   \ \ , \ T \rightarrow \infty 
\ee

          \noindent     Throughout we attribute to  $\xi(s.\chi)$ and  to $\xi \left (\frac{ 1}{2} +it  ,\chi    \right)$ and $\xi \left (t  ,\epsilon,\chi    \right)$ the same meaning. %
  As at  $\xi(t,\chi)$ and to $\xi(t,\epsilon=0,\chi)$,or at  $\xi(t)$ and to $\xi(t,\epsilon=0)$.

\noindent This article is structured as follows.
\begin{itemize}
\item In section 2  we define our angular momentum $\mathcal{L}$  connecting it with the phase variation along $t$.
\item In section 3 we give some  equivalence theorems  linking the RH or GRH   to the behavior  of $\mathcal{L}$. 

 \end{itemize}

\section{Angular Momenta}

\subsection { Phase variations and Angular Momentum of the $\xi$ function}

Derivative of the phase of $\xi(s)$ for constant $\epsilon$ with respect to $t$ :%
\begin{eqnarray}
\frac{\partial }{\partial t}  \arctan \left(  \frac {\Im \xi(s)}{ \Re \xi(s) } \right) &=&
 \frac{   1}{1+\left(  
 \frac {\Im\xi(s)}{\Re[\xi(s)]} 
  \right)^2}
   \frac{ \frac{\partial  \Im \xi(s) }{\partial t} \Re\xi(s)           -  \frac{\partial  \Re \xi(s) }{\partial t} \Im\xi(s) }{\left( \Re\xi(s) \right)^2}
\label {VarPhase} \\
&=& 
\frac{ \frac{\partial  \Im \xi(s)}{\partial t} \Re\xi(s)  -  \frac{\partial  \Re \xi(s) }{\partial t} \Im\xi(s)         }{ \left(\Re\xi(s)\right)^2 +\left( \Im\xi(s)\right)^2} = \frac{\partial \angle[\xi(s)]}{\partial t} \nonumber
\end{eqnarray}
Where $\angle[\psi ]$ means phase of $\psi$. The numerator, that determines the sign of (\ref {VarPhase}), can be seen as the angular momentum with respect to the origin of an unitary mass positioned in $(\Re\xi(s) , \Im\xi(s))$  at time $t$ for constant $\epsilon$. We will name it  $\mathcal{L} [ \xi (s) ]$, so that:
\be \label {AngMomDef}
 \mathcal{L} [ \xi (s) ] =
\det  \left(\begin{matrix} \Re\xi(s)  &  \Im\xi(s) \\ \frac{\partial}{\partial t}\Re\xi(s)  & \frac{\partial}{\partial t}\Im\xi(s)  \end{matrix}\right) 
\ee

\noindent Since $\Im \xi=0$ for $\epsilon=0$, we have $\mathcal{L}[\xi(s)]=0$ on the critical line.

\

\subsection{ First lemma on angular momentum}\label {Lemma1AngMom}
If $A(s)$ is a derivable function and $F$ is a complex constant, then:
\be \label{lemma1}
\mathcal{L}[FA(s)] =
 |F|^2 \mathcal{L}[A(s)] 
\ee

\noindent PROOF: $\angle [F  \  A ]=\angle [F ] +\angle [ A ]$ so $\frac {\partial \angle[F  \  A]}{\partial t}=\frac {\partial \angle[  A]}{\partial t}$.From (\ref {VarPhase}) (\ref{lemma1}) follows.

  \noindent            Let us apply  First Lemma %
  to \ref  {FuncEquForXiL} ( remember $  \chi_{primitive}$).
   For \cite [p.~66] {Davenport1980}  $|\tau(\chi)  |^2= q$, then  $  \left|   \frac{i^\alpha q^{\frac{1}{2}}}{\tau(\chi)}  \right| =1  $. But,  %
    at $\epsilon= 0$, for (\ref  {XiForLfunc} ), and ( \ref{FuncEquForXiL}):  %
  $$\xi(\bar{s},\bar{\chi})   e^{ - \frac{i}{2}\angle\left [  \frac{i^a q^{\frac{1}{2}}}{\tau(\chi)}  \right ]} =\xi(s,\chi)   e^{  \frac{i}{2}\angle\left [  \frac{i^a q^{\frac{1}{2}}}{\tau(\chi)}  \right ]}=\eta(t) 
    \ \  \in  \Re $$ So: %
  $\mathcal{L} [ \xi (t,\epsilon=0,\chi)]=\left | e^{ + \frac{i}{2}\angle\left [  \frac{i^a q^{\frac{1}{2}}}{\tau(\chi)}  \right ]} \right |^2 \mathcal{L}[ \eta(t )]=0$.

\subsection{ Second lemma on  angular  momentum} \label{Lemma2AngMom}
If $A(s)$ and $F(s)$ are derivable functions, then:
\be \label {lemma2}
\mathcal{L}[F(s)A(s)] =
|F(s)|^2 \mathcal{L}[A(s)] + |A(s)|^2 \mathcal{L}[F(s)] 
\ee

\noindent From  (\ref{InfiniteProdHad}):  $\xi(s) = Z^T(z)  E^T(z)= \prod_{\rho } \left(  1-\frac{z}{\rho} \right) E^T(z)$  ,where  \cite[p.~39,46]{Edwards:1974cz},   $E^T(z)=  e^{g(z)}  e^{\sum_{m=-N}^N\frac {z}{z_m} }   =\xi(0)     $, so, $\mathcal{L} \left[E^T(z)\right]=0$. Instead: %

\be  \label {Had}
\xi(z,\chi) = Z^T(z) \times E^T(z,\chi) = \left[ \prod_{m=-N'}^N \left(1- \frac {z}{z_m}   \right) \right] \times \left[   e^{g(z,\chi)}e^{ \sum_{m=-N'}^N\frac {z}{z_m} }  \right] \quad ; T\rightarrow \infty
\ee

 \noindent     Where $N$ are the zeros within $0<t<T$ while $N'$ are the zeros within $0>t>-T$.  For Dirichlet primitive characters  $g(z)=A(\chi_{primitive})+B(  \chi_{primitive}  ) z$.
    Where  $Z^T (z)=  \prod_{m=-N'}^N \left(1- \frac {z}{z_m} \right)  $, and $E^T(z,\chi)=  e^{g(z,\chi)}  e^{\sum_{m=-N'}^N\frac {z}{z_m} }        $. 
 
 \noindent %
But we state  also that
: $\mathcal{L} \left[ E^T(z,\chi)\right] =\mathcal{L} \left[ e^{g(z,\chi)}  e^{\sum_{m=-N'}^N\frac {z}{z_m} }   \right] =0 $

     \noindent  {\bf  In other words  $E^T(z,\chi)$ factor behaves like a constant   in (\ref{lemma2})  for $\xi(s,\chi)$ as well as $E^T(z)$ for   $\xi(s)$}. %

  \noindent      PROOF: from (\ref  {lemma2}), for  $\mathcal{L} [Z^T E^T(z,\chi)]$ we have (dropping :$. .  [(z,\chi)] . .$)
 \be \label {TwoDet}
 \mathcal{L} (Z^T E^T)=
 det\left( \begin{matrix}       \Re[ Z^T ]   & \Im[ Z^T ]
 \\
  \Re 
 \left[    \frac{\partial}{\partial y}          Z^T     \right]
            &  
             \Im  
\left[    \frac{\partial}{\partial y}          Z^T     \right] 
  \end{matrix}    \right) \left\{  \Re^2[  E^T ] + \Im^2[ E^T ]   \right\} + . .
  \ee
  $$
   . .+ \left\{  \Re^2[ Z^T ] + \Im^2[ Z^T]   \right\}
  det\left( \begin{matrix}       \Re[ E^T  ]   & \Im[ E^T ]
 \\
  \Re 
 \left[    \frac{\partial}{\partial y}          E^T    \right] 
            &  
             \Im
\left[    \frac{\partial}{\partial y}          E^T    \right] 
  \end{matrix}    \right) 
 $$
 
\noindent  We want to show that $det\left( \begin{matrix}       \Re[ E^T  ]   & \Im[ E^T ]
 \\
  \Re 
 \left[    \frac{\partial}{\partial y}          E^T    \right] 
            &  
             \Im
\left[    \frac{\partial}{\partial y}          E^T    \right] 
  \end{matrix}    \right) =0$   for $\xi(s,\chi)$. 
  
  \noindent Applying (\ref {BAndRoots}):

\be \label {DeEN suDeY}
 \frac{\partial}{\partial y}  E^T = \frac{\partial}{\partial y}  \left[ e^{g(z,\chi)}  e^{\sum_{m=-N'}^N\frac {z}{z_m} }  \right]= e^{g(z,\chi) + \sum_{m=-N'}^N\frac {z}{z_m} } \left(  i B(\chi) +  \sum_{m=-N}^N\frac {i}{z_m}   \right)= E^T \left(  i B(\chi) +  \sum_{m=-N'}^N\frac {i}{z_m}   \right) 
 \ee

 $$= E^T \left\{- \Im  \left[   B(\chi) +  \sum_{m=-N'}^N\frac {1}{z_m}   \right]+ i \times
\left(    \Re  \left[   B(\chi) +  \sum_{m=-N'}^N\frac {1}{z_m}   \right] \right)
\right\} =$$
$$E^T \left\{- \Im  \left[   B(\chi) +  \sum_{m=-N'}^N\frac {1}{z_m} +i \times 0   \right]\right\}=E^T \ C
$$.

\noindent  $C\in \Re$, so :%
$det\left( \begin{matrix}       \Re[ E^T  ]   & \Im[ E^T ]
 \\
  \Re 
 \left[    \frac{\partial}{\partial y}          E^T    \right] 
            &  
             \Im
\left[    \frac{\partial}{\partial y}          E^T    \right] 
  \end{matrix}    \right)=det\left( \begin{matrix}       \Re[ E^T  ]   & \Im[ E^T ]
 \\
  \Re 
 \left[    C         E^T    \right] 
            &  
             \Im
\left[   C        E^T    \right] 
 \end{matrix}    \right) = 0$.  As  $\Re[C \times E^T]=C \Re[E^T]$, and,  $\Im[C \times E^T]=C
  \Im[E^T]$.
 \noindent      %
  \noindent   So,  for both $\xi(s)$ and   $\xi(s,\chi_{primitive})$, only   determinant with $Z^T$ in (\ref {TwoDet} )  is left,  as $T \rightarrow \infty$.

\section{%
Some RH-GRH equivalences}
The infinite  products \ref{InfiniteProdHad} uses the zeros of the functions $\xi(s) \ ,  \ \xi(s,\chi)$ and the convergence in the expansion is ensured if the symmetrical  factors with respect to the real axis are joined together, that is if the zeros of the form $\rho$ and $1 - \rho$ are paired for $\xi(s)$, or are at least included symmetrically for $\xi(s,\chi)$.

\subsection { Third lemma on angular momentum}%
\label {3rdLemma} 
If all the  zeros of the $\xi$ function are  on the critical line then the sign of the angular momentum equals the sign of $\epsilon$.
That is:
\be \label {EQ1}
 \forall t : \sign(\mathcal{L} [\xi(t,\epsilon)]) = \sign(\epsilon) \quad and \  \forall t : \sign(\mathcal{L} [\xi(t,\epsilon,\chi_{primitive)}])= \sign(\epsilon)
\ee

\noindent PROOF:  %
to compute ( \ref  {lemma2} ) we can write : %
\be  \label {ZelevNm1}
\frac{\partial \xi(s)}{\partial t}=E^T(s)\frac{\partial Z^T(s)}{\partial t}
 = \xi(s) \sum_{\rho} \frac{-i}{\rho-\frac{1}{2}-\epsilon-it} 
= \xi(s) \sum_{\rho} \frac{(t-\Im[\rho)] +i(\frac{1}{2}+\epsilon-\Re[\rho])}{|\rho-s|^2} 
\ee
For $\xi(s)$  because $E^T(s)$ is constant, but, only in order to compute ( \ref  {lemma2} ), also for $\xi(s,\chi)$, even if  $E^T(s,\chi)$ is not constant (Lemma 2).
Taking now the real and imaginary parts to compute the angular momentum:
\be \label {AngMomZ}
\mathcal{L}[\xi(s)] =
 \det
\left( 
\begin{matrix}       \Re[\xi]  & \Im[\xi ] \\
  \sum_{\rho} \frac{ \Re[\xi]%
  (t-\Im\rho)  - \Im[\xi]%
   (\frac{1}{2}+\epsilon-\Re\rho)}{|\rho-s|^2} 
   &  
  \sum_{\rho} \frac{ \Re[\xi]
   (\frac{1}{2}+\epsilon-\Re\rho) + \Im[\xi ]
   (t-\Im\rho) }{|\rho-s|^2} 
\end{matrix}
\right)
\ee
so that, after cancellations:
\be \label {AngMomPolinomioN}
\mathcal{L}[\xi(s)] =|\xi(s)|^2  \sum_{\rho} \frac{  \frac{1}{2}+\epsilon-\Re[\rho]}{|\rho-s|^2}
\ee
Notice that the multiplication by $|\xi(s)|^2$ avoids divergence on each  zero.
From eq. (\ref{AngMomPolinomioN}) we see that, if $\Re\rho = \frac{1}{2}$ for each zero, %
then the  sign of  the angular momentum $\mathcal{L}[\xi(s)]$ equals the sign of $\epsilon$  everywhere.
For Lemma 2 the same applies to $\xi(s,\chi_{primitive})$.

\subsection { Fourth lemma on angular momentum%
} \label {AngMomEquTh}  

If the angular momentum of the $\xi$ function has always the same sign as $\epsilon$, then all the zeros of $\xi(s)$ are on the critical line.

\noindent PROOF: 
suppose that there is a zero $\rho^*$ with $\Re \rho^*=\frac{1}{2}+\epsilon^*$ with $\epsilon^*>0$ (by symmetry it is enough to consider this case). Eq. (\ref{AngMomPolinomioN}) becomes:

\be \label {AngMomPolinomioNbis}
\mathcal{L}[\xi(s)] =|\xi(s)|^2 \left( \frac{  \epsilon-\epsilon^*}{|\rho^*-s|^2}  + \sum_{\rho\neq \rho^*} \frac{  \frac{1}{2}+\epsilon-\Re[\rho]}{|\rho-s|^2}  \right)
\ee
Since the first term in parenthesis is the only divergent one for $s\rightarrow \rho^*$, it dominates over the rest of the sum in a neighborhood of $\rho^*$. It is then possible to choose a $\epsilon$ value such that $0<\epsilon<\epsilon^*$ where $\mathcal{L}[\xi] <0$, against the hypothesis. 
It follows that if $\sign(\mathcal{L} [\xi(s)] )= \sign(\epsilon)$ then all the zeros must be on the critical line.
For (\ref {VarPhase}) the phase variation is :

\be \label {VarFaseSeOffCLine}
  \frac {\partial \angle[\xi(t,\epsilon)]}{\partial t}  =
\frac {\mathcal{L}[\xi(s)] }{|\xi(s)|^2} 
=
\left( \frac{  \epsilon-\epsilon^*}{|\rho^*-s|^2}  + \sum_{\rho\neq \rho^*} \frac{  \frac{1}{2}+\epsilon-\Re[\rho]}{|\rho-s|^2}  \right)
 \ee
 
 For Lemma 2 the same applies to $\xi(s,\chi_{primitive})$.By Lemma 3+4, statement  (\ref {EQ1}) is equivalent to RH  and GRH.

 \subsection  { Another equivalence linked to sign of  angular momentum on whole critical strip-width.  }  \ \
 
 \noindent  Increase of $||\xi(t,\epsilon,\chi)|| 
 $, along $\epsilon$,  for $\quad\epsilon>0$, or decrease of $||\xi(t,\epsilon,\chi)|| \quad $for $\quad \epsilon<0$ are equivalent to RH for  $\xi(s,\chi_{primitive})$ and for  $\xi(s)$
because: $$   \mathcal{L}  [ \xi (t,\epsilon,\chi) ]  =   %
  \frac {1}{2}\frac{\partial ||\xi(t,\epsilon,\chi)||^2}{\partial \epsilon}$$ 
  \noindent PROOF: as $\xi(s,\chi)$  is holomorphic, due to Cauchy-Riemann relations,
$
  \forall t \quad and \quad  \forall \epsilon \quad : \quad  \mathcal{L} _\xi(t,\epsilon,\chi) = 
$

$$
   det \left(\begin{matrix} 
  \Re[\xi(  t,\epsilon,\chi   )] & \Im[\xi(  t,\epsilon,\chi )]  \\ 
  \frac{\partial \Re[\xi(  t,\epsilon,\chi )]}{\partial t} &  \frac{\partial \Im[\xi( t,\epsilon,\chi )]}{\partial t}    \end{matrix}\right) = 
  det \left(\begin{matrix} 
  \Re[\xi(t,\epsilon,\chi)] & \Im[\xi(  t,\epsilon,\chi  )]  \\ 
  -\frac{\partial \Im[\xi( t,\epsilon,\chi)]}{\partial \epsilon} &  \frac{\partial \Re[\xi( t,\epsilon,\chi )]}{\partial \epsilon}  \end{matrix}\right) = 
$$

\be \label {AngMomModQ3}
=\Re[\xi( t,\epsilon,\chi)]  \frac{\partial \Re[\xi(   t,\epsilon,\chi)]}{\partial \epsilon} +
\Im[\xi( t,\epsilon,\chi)]  \frac{\partial \Im[\xi(  t,\epsilon,\chi)]}{\partial \epsilon}  =\frac{1}{2} \frac{\partial  ||\xi( t,\epsilon,\chi )||^2}{\partial \epsilon}
\ee

  \noindent So  ( \ref  {AngMomModQ3} )  , with (\ref {VarPhase}) tel us that, if RH is true: %
  \be \label {VariazAngEModuloQ}
 \frac{\partial}{\partial t} \arctan \left(  \frac{ \Im[\xi(s.\chi) }{ \Re[\xi(s,\chi)} \right) =\frac{ \mathcal{L} [ \xi (s) ]}{||\xi(s,\chi)||^2}  =
 \frac {\frac{\partial}{ \partial \epsilon} ||\xi(s,\chi)||^2}{ 2||\xi(s,\chi)||^2} =
 \frac{1}{2} \frac{\partial}{ \partial \epsilon}  \ln(||\xi(s,\chi)||^2) 
  \ee

 \noindent As logarithm is    a monothonic function  the thesis is equivalent to Lemmas 3+4.
  And this is true  for $\xi(s)$ as well as  for $\xi(s,\chi_{primitive})$.
  Equivalence between  increasing $|| \xi(s) ||$ along $\epsilon >0$ at constant $t$ and RH is treated also in 
\cite{Lagarias1999} in wich is shown that RH is equivalent to a positivity
property of the real part of the logarithmic derivative of the $\xi(s)$ function,i.e:
$$ for \  \Re(s)>0 \rightarrow \Re\left[ \frac{\xi(s)'}{\xi(s)} \right] >0$$
\noindent Also in  \cite  {J. Sondow and C. Dumitrescu} where is shown that RH is true iff $|\xi(s)|$ is decreasing for fixed $t$ and $-\infty < \epsilon < \frac{1}{2}$, while increasing for $\frac{1}{2} < \epsilon < \infty $ . 

 \noindent The above results are proved only for RH. Here same equivalences are extended to $|| \xi(s ,\chi_{primitive}) ||$, i.e. to GRH, thanks to Angular Momentum concept  and with  Lemma 3+4, statement, i.e.   (\ref {EQ1}) .  %

\begin{figure}[!htbp]
\begin{center}
\includegraphics[width=1.0\textwidth]{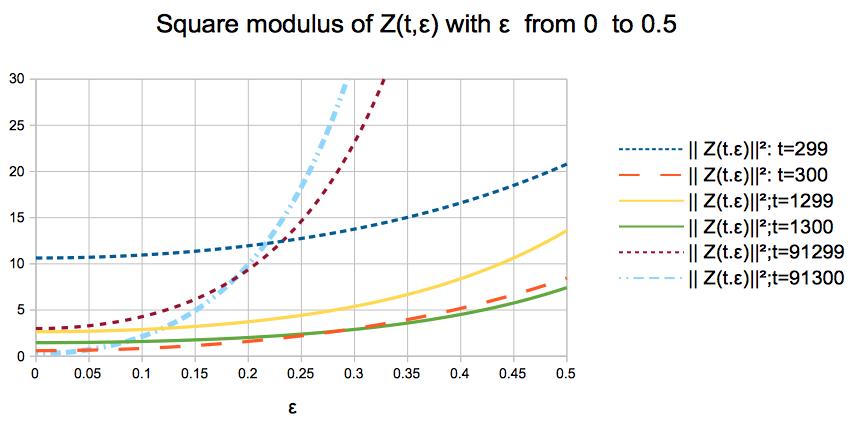} 
\caption{\small {\it    Square module of  $Z(t,\epsilon) = \frac{-\xi(t,\epsilon)}{F(t)  e^{i\epsilon\frac{\pi}{4}}} $  ( \ref  {ZSinhECosh1} )
, in the interval   $0< \epsilon < 0.5$, for different $ t$  values.
} }%
\label {Zquadro}  
\end{center}
\end{figure}

\noindent In order to see a graphical picture of this behavior for  $\xi(s)$ function, it can be computed using  {\it Wolfram Mathematica}  by RiemannSiegelZ function, or, the result is the same  ( see \cite {Giovanni Lodone 2024}),  
 by  \cite {Giovanni Lodone 2021}. See fig. \ref{Zquadro} .
 
 \noindent  Through a positive  scale factor $ F(t) =  \left(\pi/2 \right)^{0.25}t^{\frac{7}{4}} e^{- \frac{\pi}{4}  t}  $  %
 and  hyperbolic functions plus other terms here omitted.  The  expression is:
$$
\frac{-\xi(t,\epsilon)}{F(t)  e^{i\epsilon\frac{\pi}{4}}}  \approx  Z(t,\epsilon)=
2 \sum_{n=1}^N \frac{    \cosh \left[ \epsilon \  \ln\left( \sqrt{\frac {t}{2 \pi n^2}  }\right) \right]   }{\sqrt{n}}
\cos\left(   t   \ln\left(    \sqrt{\frac {t}{2  e \pi n^2} } \right) - \frac{\pi}{8}   \right) 
$$
\be + 2 i  \sum_{n=1}^N\frac{ \sinh\left[ \epsilon \  \ln\left( \sqrt{\frac {t}{2 \pi n^2}  }\right) \right]   }{\sqrt{n}}
\sin\left(   t  \  \ln\left(    \sqrt{\frac {t}{2  e \pi n^2} } \right) - \frac{\pi}{8}  \right )         +R_0(t,\epsilon) \label {ZSinhECosh1}
\ee
where $N= \left  \lfloor       \sqrt{\frac {t}{2   \pi } }   \right    \rfloor$ and: 
\be \label {R0} 
 R_0(t,\epsilon) = (-1)^{N-1} \left( \frac{2 \pi}{t}\right)^{1/4}  \left[     C_0(p)
 \right]
\ee
with $p =  \sqrt{\frac {t}{2   \pi } } -N$. In (\ref{R0}), the dependence from $\epsilon$ has been neglected  (\cite [p.~1]{Giovanni Lodone 2021}).
Besides (\cite [p.~16]{Giovanni Lodone 2021}):
\be \label {CZero}
C_0(0.5)%
\le C_0(p)= \frac{(\cos(2 \pi(p^2-p-\frac{1}{16}))}{\cos(2 \pi p)} \le \cos(\frac{\pi}{8}) 
=C_0(0)=C_0(1)
\ee

   %
\subsection  { %
A practical equivalence  based { \it ONLY} on critical line points }   \ \  \label {DerAngMom}

\noindent  Let us define ( see (\ref{FuncEquForXiL}) , and section \ref   {Lemma1AngMom}):

\be \label {EpsZero}
\eta(s,\chi) =   \left( e^{ - \frac{i}{2}\angle\left [  \frac{i^\alpha k^{\frac{1}{2}}}{\tau(\chi)}  \right ]} \xi(s,\chi) \right) ; \ for \  \epsilon=0 \quad it \quad is  \quad real \quad \forall \chi_{primitive}
\ee. 
     
     \noindent { \bf If we have :    $\left[  \frac{\partial   \mathcal{L}  [ \eta (t,\epsilon,\chi) ]  )}{\partial \epsilon} \right]_{\epsilon=0} >0$  then  RH is granted because  
  we can  extend to $\eta(t)$ 
  the %
   statement  in  %
   \cite[p.~6] {Bombieri:2000cz} }  : { \it ''The Riemann hypothesis is equivalent to the statement that all local maxima
of $  \xi(t)$ are positive and all local minima are negative, . . .'' }.
  \noindent  %
  Both $\xi(s)$ and $\eta(s,\chi_{primitive} )$ are holomorphic with real values on critical line. So we follow on with  $\eta(s,\chi_{primitive} )$. We state:

\be \label {AngMomDerivSuEps}
 \left[  \frac{\partial   \mathcal{L}  [ \eta (t,\epsilon,\chi) ]  )}{\partial \epsilon} \right]_{\epsilon=0}=
  det \left( \begin{matrix}
  \eta(t,\chi ) &- \eta'(t ,\chi)   \\ 
  \eta'(t,\chi) & -\eta''(t,\chi ) \end{matrix}\right) %
  \ee

    \noindent   PROOF:   in  \cite[p.~4]  {Giovanni Lodone Nov2024}

\noindent So, on relative maxima and minima, where $\eta'(t,\chi)=\frac{d \eta(t,\chi )}{dt} =0$,  if we  have:

  \be \label {maxminrel}
  \eta(t,\chi ) \left[   -\frac{d^2 \eta(t,\chi )}{dt^2}  \right] >0
\ee
   
\noindent       {\bf  then  RH is assured } ,  because  $\eta(t_{RELATIVE \ \ MAX},\chi_{primitive})>0$ and $\eta(t_{RELATIVE \ \ MIN}, \chi_{primitive})<0$. 
This for 
  \cite[p.~6] {Bombieri:2000cz} is equivalent to   RH.
 \noindent  All this section  applies  to $\eta(s,\chi_{primitive})$  and to $\xi(s)$  as well.


\section {Conclusions}

\noindent The  introduced  Angular Momentum / Phase derivative  of $\xi$ function concept allows to gather several  RH equivalence theorem with almost same statements  for all Dirichlet functions. In  \cite {Giovanni Lodone Nov2024}  and \cite {Giovanni Lodone Dec2024} we used the equivalence  \cite[p.~6] {Bombieri:2000cz} , here mentioned in  paragraph (3.4),  for an  attempt to prove RH for $\xi(s)$ Riemann function and for $L(s,\chi_{ odd \ primitive})$ Dirichlet functions.  %
But   Angular Momentum / Phase derivative  of $\xi$ function concept can be used also  to deduce new or known equivalences for $\xi(s)$ and to extend them easily to $\xi(s,\chi_{primitive})$. See  for example section (3.3). To our knowledge this extention  has not been published yet.

$$Acknowledgments
$$
I thank Paolo Lodone for extremely useful discussions, and, for contributing in some point of this work.
 %


\vspace{0.3cm}

\end{document}